\documentclass{article}%
\usepackage{amssymb}
\usepackage{amsmath}
\usepackage{amsfonts}
\usepackage{graphicx}
\usepackage[margin=1.25in]{geometry}%
\usepackage{ragged2e}
\providecommand{\U}[1]{\protect\rule{.1in}{.1in}}

\begin{document}

\begin{center}
{\LARGE Large time behavior of the solution to the Cauchy problem for the
discrete} {\LARGE p-Laplacian with density on infinite graphs.}

Alan A.Tedeev. Department of Mathematics, Sergo Ordzhonikidze Russian State
University for Geological Prospecting, Moscow, Russia.

e-mail: alan.tedeev2013@gmail.com, ORCID: 0009-0009-5065-4102

\end{center}
\begin{center}
\textbf{Abstract}

\begin{center} 
	\begin{minipage}{0.9\textwidth} 
		\justify 
		{\small 
			We consider the Cauchy problem for the nonstationary discrete
			p-Laplacian with inhomogeneous density $\rho(x)$ on an infinite graph which supports the Sobolev inequality.
			For nonnegative solutions when $p>2$, we prove the precise rate of
			stabilization in time, provided $\rho(x)$ is a non-power
			function. When $p>2$ and $\rho(x)$ goes to zero fast enough, we prove the universal bound. Our technique relies on suitable energy inequalities and a new embedding result.}
	\end{minipage}
\end{center}
\end{center}
Keywords$:$ Infinite Graphs $\cdot$ p-Laplacian $\cdot$ Inhomogeneous density
$\cdot$ Large time behavior $\cdot$ Universal bounds 

\noindent Mathematics Subject Classification(2020) 35R02 \textperiodcentered\ 58J35
\textperiodcentered\ 39A12

\begin{center}
\textbf{1} \textbf{Introduction}.
\end{center}

We consider nonnegative solutions to the Cauchy problem for discrete
$p$-Laplace parabolic equations with inhomogeneous density of the form

\begin{equation}
\rho(x)\frac{\partial u}{\partial t}(x,t)-\Delta_{p}u(x,t)=0\text{, }x\in
V,\text{ }t>0,\text{ }(x,t)\in S_{T}=V\times(0,T)\text{,} \tag{1.1}%
\end{equation}

\begin{equation}
u(x,0)=u_{0}(x)\geq0\text{, }x\in V\text{.} \tag{1.2}%
\end{equation}
Here $V$ is the set of vertices of the graph $G:=G(V,E,w)$ with edge set
$E\subset V\times V$ and weight $w$,
\[
\Delta_{p}u(x,t)=\frac{1}{m(x)}%
{\displaystyle\sum\limits_{y\in Y}}
\left\vert u(y,t)-u(x,t)\right\vert ^{p-2}(u(y,t)-u(x,t))w(x,y)\text{,
}N>p>1,
\]

\[
\rho(x)=\rho(d(x))\text{, }%
\]
where $d(x):=d(x,x_{0})$, $x_{0}\in V$, $x_{0}$ be fixed is combinatorial
distance in $G.$ That is, for any \\ $x$,$y$ $\in V$, the distance between them is defined as the length of the shortest path connecting $x$ and $y$. 
We use from now on for  $f:V\rightarrow$ $\mathbb{R}$ the notation
\begin{equation*}
	D_{y}f(x) = f(y) - f(x) = -	D_{x}f(y),  x,y \in V.    
\end{equation*}

We assume that the graph $G$ is simple, undirected, infinite, connected with
locally uniformly finite degree $m(x)$. The latter means that
\[
m(x):=%
{\displaystyle\sum\limits_{y\thicksim x}}
w(x,y)\leq C\text{, }%
\]
where we write $y\thicksim x$ if and only if $(x,y)\in E$. Here the weight
$w$: $V\times V\rightarrow\lbrack0,\infty)$ is symmetric, i.e.,
$w(x,y)=w(y,x)$, and strictly positive if and only if $y\thicksim x$, and
$w(x,x)=0,$ for $x\in V$. If $G$ supports the Sobolev inequality, then [2] the measure is nondegenerate i.e.

\[
\inf_{x\in V}m(x)>0\text{.}%
\]

The density function $\rho(x):V\rightarrow$ $\mathbb{R}\in\mathbb{N}$ is
positive decreasing function such that $\rho(x)\leq1$ . Additionaly, in what
follows we assume that the following assumptions hold:

$(H_{1})$: For any $s>0$ and for given $\alpha_{1}$, $\alpha_{2}$:
$0\leq\alpha_{1}\leq\alpha_{2}<$ $p$ the functions $\rho(s)s^{\alpha_{2}},$
$\rho(s)s^{\alpha_{1}}$ are increasing and decreasing respectively.

$(H_{2}):$ For any $R\geq 1$ there exists $C_{1}$ such that

\[
\mu (B(R)):=\sum\limits_{x\in B(R)}m(x)\leq C_{1}R^{N}.
\]%
The purpose of the paper is to get precise temporal decay estimates of the
solution to (1.1), (1.2) under the assumptions on the geometry of graph and
behavior of $\rho(x)$.

Define%

\[
\left\Vert f\right\Vert _{l^{q}(U)}^{q}=%
{\displaystyle\sum\limits_{x\in U}}
\left\vert f(x)\right\vert ^{q}m(x)\text{, }\left\Vert f\right\Vert
_{l^{\infty}(U)}=\sup_{x\in U}\left\vert f(x)\right\vert \text{, }\mu(U)=%
{\displaystyle\sum\limits_{x\in U}}
m(x)\text{,}%
\]
where all infinite sums are absolutely convergent.

\textbf{Definition 1.1} \textit{We say that }$G$\textit{\ satisfies the
Sobolev inequality for some given }$p\geq1$, $N>p$ \textit{if for any }$v>0$
\textit{and any finite subset} $U\subset V$\textit{ with} $\mu%
(U)=v$\textit{ there exists constant }$C$ \textit{independent of} $v$
\textit{such that} \textit{the following inequality holds true}%
\begin{equation}
\left(
{\displaystyle\sum\limits_{x\in U}}
\left\vert f(x)\right\vert ^{p^{\ast}}m(x)\right)  ^{1/p^{\ast}}\leq
C\left(
{\displaystyle\sum\limits_{x,y\in U}}
\left\vert f(y)-f(x)\right\vert ^{p}w(x,y)\right)  ^{1/p}, \tag{1.3}%
\end{equation}
\textit{for all real-valued functions} $f:V\rightarrow \mathbb{R}$ \textit{such that} $f(x)=0$ if $x\notin
U$. \textit{Here }$p^{\ast}=Np/(N-p)$.	

\textbf{Remark 1.1} Here and hereafter we write
\[
\left\Vert u(t)\right\Vert _{\infty}:=\left\Vert u(x,t)\right\Vert
_{l^{\infty}(V)}\text{,}%
\]
by $C$ we denote a generic constant which depends only on parameters
$p,\alpha_{1},\alpha_{2},N$ of the problem and may vary from line to line.

\textbf{Remark 1.2} Note that (1.3) follows from the $N-$isoperimetric
inequality of $G$ (see, for example, [9], [17], [34]): for any $v>0$ and any
finite subset $U\subset V$ with $\mu(U)=v<\infty$ there exists constant
$C$ independent of $v$ such that%

\begin{equation}
\mu(U)^{(N-1)/N}\leq C\mu(\partial U)\text{,} \tag{1.4}%
\end{equation}
where $\partial U$ is the set of all edges with one endpoint in $U$ and another
one in $V\backslash U$, \\ $\mu(\partial U):=%
{\displaystyle\sum\limits_{x,y\in U}}
w(x,y).$
Note also that (1.4) implies the bound from below on the volume%

\[
\mu(B(R)):=%
{\displaystyle\sum\limits_{x\in B(R)}}
\rho(d(x))m(x)\geq CR^{N}\text{,}%
\]
where $B(R)$ means a ball of radius $R$: $\left\{  x\in V:\text{ }d(x)\leq
R\right\}  $. Then, assuming $(H_{2})$ it follows that
$\mu(B(R))\asymp R^{N}$.

The class of graphs satisfying this volume growth condition is broad (see, e.g., the survey [39] for examples). For concreteness, we will focus on the canonical case of $G = \mathbb{Z}^N$ in our analysis (see Example 1.1 below).

\textbf{Remark 1.3} By applying the H\"{o}lder inequality from (1.3) it is
easy to obtain the following Faber-Krahn inequality%

\begin{equation}%
{\displaystyle\sum\limits_{x\in U}}
\left\vert f(x)\right\vert ^{p}m(x)\leq C\mu(U)^{p}%
{\displaystyle\sum\limits_{x,y\in U}}
\left\vert f(y)-f(x)\right\vert ^{p}w(x,y). \tag{1.5}%
\end{equation}
Indeed,
\[%
{\displaystyle\sum\limits_{x\in U}}
\left\vert f(x)\right\vert ^{p}m(x)\leq\left(
{\displaystyle\sum\limits_{x\in U}}
\left\vert f(x)\right\vert ^{p^{\ast}}m(x)\right)  ^{p/p^{\ast}}\mu(U)^{p}%
\]

\[
\leq C\mu(U)^{p}%
{\displaystyle\sum\limits_{x,y\in U}}
\left\vert f(y)-f(x)\right\vert ^{p}w(x,y)\text{.}%
\]
Thus, (1.5) occurs.

\textbf{Definition 1.2} \textit{We say that }$u\in L^{\infty}(0,T;l^{r}%
(V))$\textit{ for some }$r\geq1$\textit{ and any }$T>0$ \textit{is a solution
to (1.1), (1,2) if }$u(x,t)\in C^{1}([0,T])$\textit{ for every $x\in V$ and (1.1), (1.2) are satisfied in classical point wise sense}.

Note that a solution $u$ of the problem (1.1), (1.2) exists. For this we refer the reader to [33] and [2].

In what follows, when a function is defined only on the integers, we extend it to the positive real axis by linear interpolation. For brevity, we use the same notation for both the original function and its extension. To formulate our first result we denote

\begin{equation}
	\begin{split}
		\omega(S) &= S^{p}\rho(S)\text{, } \\
		\Psi(S) &= \omega(S)S^{(N-p)(p-1)/p}\text{, } \\
		\Phi(S) &= \left( S\omega(\Psi^{(-1)}(S^{-1/p})) \right)^{1/(p-1)}\text{.}
	\end{split}
	\tag{1.6}%
\end{equation}

The functions \(\omega\), \(\Psi\), and \(\Phi\) defined in (1.6) combine the density \(\rho(S)\), the volume growth \(R^N\), and the parameter of nonlinearity. They are introduced to establish a weighted Gagliardo-Nirenberg inequality (see Lemma 2.2), which is crucial to prove Theorem 1.1.

\textbf{Theorem 1.1} \textit{Let }$u(x,t)$\textit{ be the nonnegative
solution of the problem (1.1), (1.2) in }$S_{T}=V\times(0,T)$\textit{ for any
}$T>0$\textit{ and }$G$\textit{ supports the Sobolev inequality (1.3) and
}$(H_{1})$, $(H_{2})$ \textit{hold}. \textit{ Assume that }$u_{0}\rho\in
l^{1}(V),$ $u_{0}\in l^{\infty}(V)$\textit{ and that the inverse }$\Phi
^{(-1)}$ to $\Phi$ \textit{is convex. Then for any }$t>t_{0}$, \textit{where
}$t_{0}=t_{0}(u_{0})$\textit{ is large enough, the following estimate holds
true}%

\begin{equation}
\left(  \omega^{(-1)}(\left\Vert u(t)\right\Vert _{\infty}^{p-2}t)\right)
^{N-p}\left\Vert u(t)\right\Vert _{\infty}^{p-1}tM(0)\leq C\text{,} \tag{1.7}%
\end{equation}
\textit{where}

\begin{equation*}
	M(0) := \sum_{x \in V} \rho(x) u_0(x) \, m(x).
\end{equation*}
For the continuous analog of Theorem 1.1 see [14]. However, as we will see the
proof of (1.7) is more delicate as a matter of discrete character of the problem.

\textbf{Example 1.1}  Let $G=\mathbb{Z}^{N}$, $\rho(x)=\left(  \left\vert
x_{1}\right\vert +..+\left\vert x_{N}\right\vert \right)  ^{-\alpha}$, where
$0<\alpha<p<N$, then for any $t>t_{0}$ from (1.7) it follows that%

\begin{equation}
\left\Vert u(t)\right\Vert _{\infty}\leq CM(0)^{\left(  p-\alpha\right)
/H}t^{-\left(  N-\alpha\right)  /H}\text{,} \tag{1.8}%
\end{equation}
where%

\begin{equation*}
	H:=(N-\alpha)(p-2)+p-\alpha=\lambda-\alpha(p-1),    \lambda=N(p-2)+p
\end{equation*}

\textbf{Example 1.2} Let $\rho(s):=s^{-\alpha}\left(  \log s\right)
^{\beta}$. Then $\tau=\omega(s)=s^{p-\alpha}\left(  \log s\right)  ^{\beta}$,
$0\leq\alpha<p$, $\beta>1$.

We have $s=\omega^{(-1)}(\tau)=(p-\alpha)^{\beta/(p-\alpha)}\tau
^{1/(p-\alpha)}\left(  \log\tau\right)  ^{-\beta/(p-\alpha)}(1+\epsilon
(\tau))$. 
\noindent Indeed, 

\begin{align*}
	\tau &= (p-\alpha)^{\beta} \tau \left( \log s \right)^{-\beta} (1+\epsilon(\tau)) \\
	&\quad \times \left[ \frac{\beta}{p-\alpha} \log\left(p-\alpha\right) 
	+ \frac{1}{p-\alpha} \log\tau 
	- \frac{\beta}{p-\alpha} \log\log\tau 
	+ \log(1+\epsilon(\tau)) \right]^{\beta} \\
	&\quad \times \left[ 1 + \epsilon(\omega^{(-1)}(\tau)) \right].
\end{align*}
Since, $\tau\rightarrow\infty$ as $s\rightarrow\infty$, it follows that
$\omega^{(-1)}(\tau)\rightarrow\infty$ for $\tau\rightarrow\infty$. Thus%

\[
(1+\epsilon(\tau))^{\frac{p-\alpha}{\beta}}\left[  1+o(1)+(p-\alpha)\left(
\log\tau\right)  ^{-1}\log(1+\epsilon(\tau))\right]  \rightarrow1\text{ as
}\tau\rightarrow\infty\text{.}%
\]
So, (1.7) for large enough $t$ reads as follows%

\begin{equation*}
\left\Vert u(t)\right\Vert _{\infty}\left(  \log(t\left\Vert u(t)\right\Vert
_{\infty}^{p-2})\right)  ^{-\beta(N-p)/H}\leq CM(0)^{(p-\alpha)/H}t^{-\left(
N-\alpha\right)  /H}\text{.} 
\end{equation*}
Thus, we have by maximum principle

\begin{equation}
	\left\Vert u(t)\right\Vert _{\infty }\leq C\left( \log (t\left\Vert
	u_{0}\right\Vert _{\infty }^{p-2})\right) ^{\beta (N-p)/H}M(0)^{(p-\alpha
		)/H}t^{-(N-\alpha )/H}.  \tag{1.9}
\end{equation}

\textbf{Remark 1.4} Note that $H=0$ in (1.8), (1.9) corresponds to the
critical $\alpha^{\ast}=\lambda/(p-1)$. In the continuous case for \ $\alpha>$
$\alpha^{\ast}$ the interface blow-up occurs [40]. However, for the discrete
$p$-Laplace equation the finite speed of propagation does not hold (see
comments in [33], [2]). So, in the discrete setting there is no analog of such sort of result.

\textbf{Theorem 1.2} \textit{ Let }$u(x,t)$\textit{\ be a nonnegative solution
of (1.1) in }$S_{T}$\textit{\ for any }$T>0$. \textit{Assume that
(1.3) holds and $p>2$. Then provided}

\begin{equation}%
{\displaystyle\sum\limits_{x\in V}}
\rho(x)^{N(p-1+\nu)/\left(  \lambda+p\nu\right)  }m(x)<\infty,\text{ }%
\nu>0\text{ \textit{is large enough,}} \tag{1.10}%
\end{equation}
\textit{the following universal bound holds true for all }$t>0$%

\begin{equation}
\left\Vert u(t)\right\Vert _{\infty}\mathit{\leq Ct}^{-1/(p-2)}%
\text{\textit{,}} \tag{1.11}%
\end{equation}
\textit{where }$C$\textit{ is independent of }$u_{0}(x)$.

\textbf{Remark 1.5} Let us prove that if  $\alpha >p$ and $\rho (s)s^{\alpha}$  is nonincreasing for $s>1$, then (1.10) holds. 

\textit{Proof.} The argument mirrors that of Lemma 2.1 (see Section 2). We have for $n$ large
enough

\[
\sum\limits_{x\in B(n)}\rho (d(x))^{N(p-1+\nu )/(\lambda +p\nu )}m(x)\leq
C\sum\limits_{k=1}^{n}\left( \rho (k)k^{\alpha }k^{-\alpha }\right)
^{N(p-1+\nu )/(\lambda +p\nu )-1}k^{N}
\]

\[
\leq C\sum\limits_{k=1}^{n}k^{-\alpha (N(p-1+\nu )/(\lambda +p\nu ))-1}k^{N}%
\text{.}
\]

Let $\beta =\alpha (N(p-1+\nu )/(\lambda +p\nu ))$. Then $\beta >N$ provided
for large enough $\nu $. Indeed, assume $\alpha =p+\varepsilon $. Then,
choosing $\nu $ as follows 

\[
\nu >\frac{(N-p)(p-2)}{\varepsilon }-p+1\text{,}
\]
we conclude that $\beta >N$. The latter means that

\[
\sum\limits_{x\in B(n)}\rho (d(x))^{N(p-1+\nu )/(\lambda +p\nu )}m(x)\leq
\sum\limits_{k=1}^{\infty }k^{-\beta -1+N}<\infty \text{,}
\]
where $C$ is independent of $n$. Letting $n\rightarrow \infty $, we arrive
to the desired result.

Our approach relies on the energy method of the O.A. Ladyzhenskaya, N.N.
Uraltseva, E. DiBenedetto in the form proposed in [4, 5], which is flexible
enough to apply it in the discrete setting as well (see [2] by Andreucci and Tedeev). As an important ingredient in the proof of Theorem 1.1 we prove (see Lemma 2.2 and 2.4) the discrete version of the Gagliardo-Nirenberg inequality for non-power $\rho$ by using the Sobolev inequality. For the weighted Gagliardo-Nirenberg inequality
we also refer to [32]. The continuous version of that inequality was proven in [14], [3] by using among others the Hardy inequality. Since the Sobolev
inequality does not imply the Hardy inequality, our proof is of independent
interest even in continuous case.

In continuous setting of (1.1), (1.2) the universal bound time decay estimate phenomena was proven in [40]. The surprising properties related
to the qualitative behavior of degenerate parabolic equations in the
continuous setting goes back to pioneering papers by Kamin and Rosenau [23,24]. For the further results in that direction we refer the reader to [19], [25,29] and references therein. The universal bound phenomena was investigated in several
papers (see, for example, [3], [5] and references therein).

To the best of our knowledge, the universal bound in weighted case for degenerate parabolic equation on infinite graphs never were
treated before. However, related problems have been studied in multiple settings.

Several results have been obtained for the stationary (elliptic) case of the $p$-Laplacian on infinite graphs, where existence, qualitative properties, and functional inequalities have been studied in various contexts. For instance, Holopainen and Soardi [20, 21] established foundational results for $p$-harmonic functions on graphs and manifolds, including Liouville-type theorems. Keller and Lenz [26] investigated unbounded Laplacians on graphs, focusing on their spectral properties. More recently, Biagi, Meglioli, and Punzo [6] proved a Liouville theorem for elliptic equations with potentials on infinite graphs under volume growth conditions. For $p$-Laplacian elliptic inequalities, see the work of Ting and Feng [41].

In a closely related direction, functional inequalities and embedding theorems on infinite graphs have been investigated, as they are crucial in deriving qualitative properties of solutions. Notably, Pinchover and Tintarev [35] studied positive solutions of $p$-Laplacian-type equations, while Saloff-Coste [38] provided a comprehensive treatment of Sobolev-type inequalities in discrete settings. Mourrat, Otto [32] developed Nash inequalities and heat kernel bounds for degenerate environments, which are closely related to embedding results.

For the blow-up phenomena in finite graphs we refer the reader to [11]. Note that for $p>2$ the mass conservation property for solutions was established in [15]. In the linear case (i.e., when $p = 2$), the qualitative theory becomes significantly more involved. For foundational results on analysis of weighted graphs, including heat kernel estimates and geometric conditions for stochastic completeness, we refer to Grigor'yan [17]. Andres, Deuschel, and Slowik [1] derived heat kernel estimates for random walks with degenerate weights. Earlier works by Chung, Grigor'yan, and Yau [9] connected higher eigenvalues to isoperimetric inequalities on graphs and manifolds. 

The Cauchy problem for the heat equation on infinite graphs—whose geometric structures are characterized by Faber-Krahn or isoperimetric inequalities—has been studied in [1], [9], [11]. For stationary $p$-Laplacian on discrete graphs see [37]. 

For completeness, we mention that the \textit{p}-Laplacian on finite graphs has been studied in applications such as image processing (Elmoataz et al. [15]) and spectral theory (Cardoso and Pinheiro [8]). The dynamics of solutions, including extinction and positivity, were investigated by Lee and Chung [30] and Xin et al. [42].

The structure of the paper is as follows: In Section 2 we give auxiliary
results that we use in proofs of Theorems 1.1-1.2. Sections 3,4 are devoted
to the proofs of Theorems 1.1-1.2 respectively.

\begin{center}
\textbf{2} \textbf{Preliminary results.}
\end{center}

\bigskip\ Denote%

\[
E_{q}:=%
{\displaystyle\sum\limits_{x\in V}}
\left\vert f\right\vert ^{q}\rho(x)m(x),\text{ }D_{p}:=%
{\displaystyle\sum\limits_{x,y\in V}}
\left\vert D_{y}f(x)\right\vert ^{p}w(x,y)\text{,}%
\]

\[
\text{ }\psi(S)=\omega(S)S^{(N-p)(p-r)/p}\text{, }\varphi(S):=\left(
S\omega\left(  \psi^{(-1)}(S^{-r/p})\right)  \right)  ^{(q-r)/(p-r)}\text{,
\ }S>0.
\]

\textbf{Lemma 2.1} \textit{Assume that $(H_{2})$ holds. Then for any $0<\beta <N$ and $n > n_{0},$ where $n_{0}$ is large enough}

\begin{equation}
	\sum\limits_{x\in B(n)}\mathit{d(x)}^{-\beta }\mathit{m(x)\leq Cn}^{N-\beta
	}\text{,} \tag{2.1}
\end{equation}
\textit{and for any} $\beta >N$ \textit{and} $n > n_{0},$

\[
\sum\limits_{x\in V}d(x)^{-\beta }m(x)<\infty. \tag{2.2}
\]

\textit{Proof}. Let $0<\beta <N$. We have with $R=n$ that

\[
\sum\limits_{x\in B(n)\diagdown B(2)}d(x)^{-\beta
}m(x)=\sum\limits_{k=2}^{n}\sum\limits_{k\leq d(x)\leq k+1}d(x)^{-\beta
}m(x)
\]

\[
\leq \sum\limits_{k=2}^{n}k^{-\beta }\sum\limits_{k\leq d(x)\leq k+1}m(x)%
\text{.}
\]%
By Abel's summation by parts formula we have

\[
\sum\limits_{k=j}^{n}r_{k}b_{k}=A_{n}b_{n}-A_{j-1}b_{j}+\sum%
\limits_{k=j}^{n}A_{k}(b_{k}-b_{k+1})\text{,}
\]%
where

\[
A_{n}=\sum\limits_{k=1}^{n}r_{k}\text{.}
\]%
Choose now $j=2$ and set

\[
r_{k}=\sum\limits_{k\leq d(x)\leq k+1}m(x)\text{, }b_{k}=k^{-\beta }\text{.}
\]%
Then

\[
A_{n}b_{n}-A_{1}b_{2}=\left[ \mu (n+1)-\mu (1)\right] n^{-\beta }\leq
C_{1}(n+1)^{N-\beta }\text{.}
\]%
Thus, 

\[
\sum\limits_{k=2}^{n}k^{-\beta }\sum\limits_{k\leq d(x)\leq k+1}m(x)\leq
C_{1}(n+1)^{N-\beta }
\]

\[
+\sum\limits_{k=2}^{n}(k^{-\beta }-(k+1)^{-\beta })\sum\limits_{1\leq
	d(x)\leq k+1}m(x)
\]

\[
\leq C_{1}(n+1)^{N-\beta }+C\sum\limits_{k=1}^{n}k^{-\beta -1+N}\leq
Cn^{N-\beta }\text{.}
\]%
Finally,

\[
\sum\limits_{x\in B(n)}d(x)^{-\beta }m(x)=\sum\limits_{x\in
	B(2)}d(x)^{-\beta }m(x)+\sum\limits_{x\in B(n)\diagdown B(2)}d(x)^{-\beta
}m(x)
\]

\[
\leq C+Cn^{N-\beta }\leq Cn^{N-\beta }
\]
for $n>n_{0}$ large enough. This proves (2.1).  Working analogously, we have
for $\beta >N$

\[
\sum\limits_{x\in B(n)}d(x)^{-\beta }m(x)\leq C\sum\limits_{k=1}^{\infty
}k^{-\beta -1+N}\leq C\text{.}
\]
Letting $n\rightarrow \infty $, we arrive at the desired result. 

Lemma 2.1 is proved.  $\square$

\textbf{Lemma 2.2} \textit{Assume that }$G$ \textit{supports the Sobolev
inequality} \textit{assumptions on }$\rho$: $(H_{1})$ \textit{and} $(H_{2}%
)$\textit{ hold. Then for }$0<r<q<p$\textit{ the following Gagliardo-Nirenberg inequality holds true }%

\begin{equation}
\mathit{E}_{q}\mathit{\leq C\left(  D_{p}^{q-r}\omega\left(  \psi
^{(-1)}\left(  D_{p}^{-r/p}E_{r}\right)  \right)  ^{q-r}E_{r}^{p-q}\right)
^{1/(p-r)},} \tag{2.3}%
\end{equation}
provided%

\begin{equation}
\psi^{(-1)}\left(  D_{p}^{-r/p}E_{r}\right)  \geq1. \tag{2.4}%
\end{equation}

\textit{Proof}. We have

\begin{equation}
E_{q}=%
{\displaystyle\sum\limits_{x\in B(R)}}
\left\vert f\right\vert ^{q}\rho(x)m(x)+%
{\displaystyle\sum\limits_{x\in V\diagdown B(R)}}
\left\vert f\right\vert ^{q}\rho(x)m(x)=I_{1}(R)+I_{2}(R)\text{,}
\tag{2.5}
\end{equation}
where $R$ be an integer that will be selected below. By applying H\"{o}lder's and Sobolev's inequalities we get

\begin{equation}
I_{1}(R)\leq\left(
{\displaystyle\sum\limits_{x\in B(R)}}
\left\vert f\right\vert ^{p^{\ast}}m(x)\right)  ^{q/p\ast}\left(
{\displaystyle\sum\limits_{x\in B(R)}}
\rho(x)^{p^{\ast}/(p^{\ast}-q)}m(x)\right)  ^{(p^{\ast}-q)/p^{\ast}}.
\tag{2.6}%
\end{equation}
The first term in the right-hand side of (2.6) we bound above by the Sobolev inequality%

\[
\left(
{\displaystyle\sum\limits_{x\in B(R)}}
\left\vert f\right\vert ^{p^{\ast}}m(x)\right)  ^{q/p\ast}\leq\left(
{\displaystyle\sum\limits_{x\in V}}
\left\vert f\right\vert ^{p^{\ast}}m(x)\right)  ^{q/p\ast}%
\]

\begin{equation}
\leq C\left(
{\displaystyle\sum\limits_{x,y\in V}}
\left\vert D_{y}f(x)\right\vert ^{p}w(x,y)\right)  ^{q/p}:=CD_{p}^{q/p}.
\tag{2.7}%
\end{equation}
For the second term we have

\[
\left(
{\displaystyle\sum\limits_{x\in B_{R}}}
\rho(x)^{p^{\ast}/(p^{\ast}-q)}m(x)\right)  ^{(p^{\ast}-q)/p^{\ast}%
}=\left(
{\displaystyle\sum\limits_{x\in B_{R}}}
\rho(x))^{Np/h(q)}m(x)\right)  ^{h(q)/Np}%
\]%
\[
=\left(
{\displaystyle\sum\limits_{x\in B(R)}}
\left[  \rho(x)\left(  d(x)\right)  ^{p}\right]  ^{Np/h(q)}d(x)^{-Np^{2}%
/h(q)}m(x)\right)  ^{h(q)/Np}%
\]%
\begin{equation}
\leq\omega(R)\left(
{\displaystyle\sum\limits_{x\in B(R)}}
d(x)^{-Np^{2}/h(q)}m(x)\right)  ^{h(q)/Np}\text{, } \tag{2.8}%
\end{equation}
where $h(q):=N(p-q)+qp$. Here we have used an assumption $\left(
H_{1}\right)  $. Since $q<p$ we have

\[
N>\frac{Np^{2}}{h(q)}\text{ :}=\beta\text{.}%
\]
Thus, using assumption $(H_{2})$, we get from (2.6)-(2.8):

\begin{equation}
I_{1}(R)\leq CD_{p}^{q/p}\rho(R)R^{h(q)/p}. \tag{2.9}%
\end{equation}
To estimate $I_{2}$ we apply H\"{o}lder's and Sobolev's inequalities to obtain%

\[
I_{2}(R)\leq\left(
{\displaystyle\sum\limits_{x\in V}}
\left\vert f\right\vert ^{p^{\ast}}m(x)\right)  ^{(q-r)/(p^{\ast}-r)}%
\]

\[
\times\left(
{\displaystyle\sum\limits_{\left\vert dx\right\vert >R}}
\left\vert f\right\vert ^{r}\rho\rho^{(q-r)/(p^{\ast}-q)}m(x)\right)
^{(p^{\ast}-q)/(p^{\ast}-r)}%
\]

\begin{equation}
\leq CD_{p}^{N(q-r)/h(r)}E_{r}^{h(q)/h(r)}\rho^{(q-r)(N-p)/h(r)}(R)\text{.}
\tag{2.10}%
\end{equation}
Let us choose $R$ as follows

\begin{equation}
D_{p}^{q/p}\rho(R)R^{h(q)/p}\geq D_{p}^{N(q-r)/h(r)}E_{r}^{h(q)/h(r)}%
\rho^{(q-r)(N-p)/h(r)}(R)\text{. } \tag{2.11}%
\end{equation}
Observing that%

\[
1-\frac{(q-r)(N-p)}{N(p-r)+pr}=\frac{h(q)}{h(r)},
\]

\[
\frac{N(q-r)}{h(r)}-\frac{q}{p}=\frac{rh(q)}{ph(r)}=\frac{rh(q)}%
{ph(r)}\text{,}%
\]
(2.11) reads

\[
\psi(R)=R^{h(r)/p}\rho(R)\geq D_{p}^{-r/p}E_{r}\text{.}%
\]
Thus, for

\begin{equation}
R\geq\psi^{(-1)}(D_{p}^{-r/p}E_{r}).\tag{2.12}%
\end{equation}
The second term in (2.11) is bounded above by the first term. In the same time
the term $D_{p}^{q/p}\rho(R)R^{h(q)/p}$ in the left-hand side of (2.11) is
bounded above by $\left(  D_{p}^{q-r}\omega(R)^{q-r}E_{r}^{p-q}\right)
^{1/(p-r)}$. Indeed,%

\[
D_{p}^{q/p}\rho(R)R^{h(q)/p}=\left(  D_{p}^{q-r}\omega(R)^{q-r}E_{r}%
^{p-q}\right)  ^{1/(p-r)}%
\]

\[
\times\left(  D^{r/p}E_{r}^{-1}\left(  \psi(R)\right)  ^{-1}\right)
^{(p-q)/(p-r)}\leq\left(  D_{p}^{q-r}\omega(R)^{q-r}E_{r}^{p-q}\right)
^{1/(p-r)}\text{.}%
\]
Thus, collecting (2.8)-(2.11), we arrive at

\begin{equation}
E_{q}\leq C\left(  D_{p}^{q-r}\omega(R)^{q-r}E_{r}^{p-q}\right)
^{1/(p-r)}\text{.} \tag{2.13}%
\end{equation}
Finally, by (2.4), we may choose $R=n$ as follows%

\[
n=\left[  \psi^{(-1)}(D_{p}^{-r/p}E_{r})\right]
\]
to get the desired result. Here $\left[  r\right]  $ is an integer part of
$r.$

Lemma 2.2 is proved. $\square$

\textbf{Lemma 2.3}. \textit{Assume that the function }$\varphi^{(-1)}(S)$
\textit{is convex. Then for some given positive constant }$K$ \textit{under
the conditions of Lemma 2.2 \ we have }
\begin{equation}
\mathit{KE}_{q}\mathit{\leq\varepsilon D}_{p}\mathit{+C}_{\varepsilon
}K\mathit{E}_{r}^{q/r}\mathit{\varphi}_{1}^{(-1)}\left(  \frac{K}%
{E_{r}^{(p-q)/r}}\right)  \text{, }\varphi_{1}(S)=\frac{\varphi^{(-1)}(S)}%
{S}.\tag{2.12}%
\end{equation}
\textit{Proof}. From Lemma 2.2 it follows that%

\begin{equation}
KE_{q}\leq CKE_{r}^{(q-r)/r}E_{r}^{p/r}\varphi\left(  D_{p}^{-r/p}%
E_{r}\right)  . \tag{2.13}%
\end{equation}
Denote%

\[
S=D_{p}E_{r}^{-p/r}\text{. }%
\]
Then, using convexity of $\varphi^{(-1)}$ we have for any positive constants
$A$, $B$, $\varepsilon:$%

\begin{equation}
\frac{\varphi^{(-1)}(A)}{A}B\leq\varepsilon\varphi^{(-1)}(B)+C_{\varepsilon
}\varphi^{(-1)}(A).\tag{2.14}%
\end{equation}
Choose in (2.14):%

\[
\varphi_{1}(A)=CE_{r}^{(q-r)/r}K\text{, }B=\varphi\left(  \frac{D_{p}}%
{E_{r}^{p/r}}\right)
\]
to get%

\[
KE_{q}\leq\varepsilon D_{p}+C_{\varepsilon}KE_{r}^{q/r}\left(  \varphi
_{1}^{(-1)}\left(  \frac{K}{E_{r}^{(p-q)/r}}\right)  \right)  \text{.}%
\]
Lemma 2.2 is proved. $\square$

\textbf{Example 2.1}. Let $\rho(R)=R^{-\alpha},$ $0\leq\alpha<p$. Then under
the condition (2.2) we have %

\[
\omega(R)=R^{p-\alpha},\psi(R)=R^{(h(r)-p\alpha)/p},
\]

\[
\left(  D_{p}^{q-r}\omega\left(  \psi^{(-1)}\left(  D_{p}^{-r/p}E_{r}\right)
\right)  ^{q-r}E_{r}^{p-q}\right)  ^{1/(p-r)}=
\]

\[
=\left(  E_{r}^{(h(q)-p\alpha)}D_{p}^{(N-\alpha)(q-r)}\right)
^{1/(h(r)-p\alpha)}.
\]
Therefore,%

\begin{equation}
E_{q}\leq C\left(  E_{r}^{(h(q)-p\alpha)}D_{p}^{(N-\alpha)(q-r)}\right)
^{1/(h(r)-p\alpha)}. \tag{2.15}%
\end{equation}
In particular, if $\alpha=0$, then (2.15) is the classical Gagliardo-Nirenberg inequality%

\begin{equation}
E_{q}\leq C\left(  E_{r}^{h(q)}D_{p}^{N(q-r)}\right)  ^{1/h(r)}. \tag{2.16}%
\end{equation}
\textbf{Remark 2.1.} In what follows in the proof of Theorem 1.1. we will use
special choice of parameters $q$ and $r$:%

\[
q=2\text{, }r=1\text{.}%
\]
Then functions $\psi$ and $\varphi$ are transformed into the form%

\[
\psi(R)\rightarrow\Psi(R)=\omega(R)R^{(N-p)(p-1)/p}\text{, }%
\]

\[
\varphi(R)\rightarrow\Phi(R)=\left(  R\omega\left(  \left(  \Psi
^{(-1)}(R^{-1/p}\right)  \right)  \right)  ^{1/(p-1)}\text{.}%
\]
Since by our assumption $\Phi^{(-1)}(R)$ is a convex function by Lemmas 2.1
and 2.2 we have%

\begin{equation}
KE_{2}\leq\varepsilon D_{p}+C_{\varepsilon}KE_{1}^{2}\left(  \Phi_{1}%
^{(-1)}\left(  \frac{K}{E_{1}^{p-2}}\right)  \right)  \text{, }\Phi
_{1}(S):=\frac{\Phi^{(-1)}(S)}{S}\text{.} \tag{2.17}%
\end{equation}

\textbf{Lemma 2.3.} \textit{Let }$u(x,t)$\textit{ be the solution of (1.1),
(1.2) in }$S_{\infty}$\textit{. Then, if }$\left\Vert \rho u_{0}\right\Vert
_{1}<\infty,$ $\left\Vert u_{0}\right\Vert _{\infty}<\infty$\textit{, for any
}$t>0$%
\begin{equation}%
{\displaystyle\sum\limits_{x,y\in V}}
\left\vert D_{y}u(x,t)\right\vert ^{p}\mathit{w(x,y)\leq}\left(
\mathit{2/t}\right)  \mathit{\left\Vert u_{0}\right\Vert _{\infty}}%
{\displaystyle\sum\limits_{x\in V}}
\mathit{\rho u}(x,t)\mathit{d}_{w}\mathit{(x).}\tag{2.18}%
\end{equation}

\textit{Proof}. Multiplying both sides of the equation by $u_{t}$ and
integrating by parts in time, we have%

\begin{equation}%
{\displaystyle\sum\limits_{x\in V}}
\rho u_{t}^{2}(x,t)m(x)=-\frac{1}{2p}\frac{d}{dt}%
{\displaystyle\sum\limits_{x,y\in V}}
\left\vert D_{y}u(x,t)\right\vert ^{p}w(x,y). \tag{2.19}%
\end{equation}
From (2.19) it follows that for any $0<t_{1}<t_{2}$%

\begin{equation}%
{\displaystyle\sum\limits_{x,y\in V}}
\left\vert D_{y}u(x,t_{2})\right\vert ^{p}w(x,y)\leq%
{\displaystyle\sum\limits_{x,y\in V}}
\left\vert D_{y}u(x,t_{1})\right\vert ^{p}w(x,y)\text{.} \tag{2.20}%
\end{equation}
Next, multiply both sides of (1.1) by $u$ and integrating by parts in time in the time interval $(t/2, t)$ to get%

\[%
{\displaystyle\sum\limits_{x\in V}}
\rho u^{2}(x,t)m(x)=%
{\displaystyle\sum\limits_{x\in V}}
\rho u^{2}(x,t/2)m(x)-%
{\displaystyle\int\limits_{t/2}^{t}}
{\displaystyle\sum\limits_{x,y\in V}}
\left\vert D_{y}u(x,t)\right\vert ^{p}w(x,y)dt\text{.}%
\]
Then, due to (2.20)%

\[%
{\displaystyle\sum\limits_{x,y\in V}}
\left\vert D_{y}u(x,t)\right\vert ^{p}w(x,y)\leq2/t%
{\displaystyle\sum\limits_{x\in V}}
\rho u^{2}(x,t/2)m(x)
\]

\[
\leq2/t\left\Vert u_{0}\right\Vert _{\infty}%
{\displaystyle\sum\limits_{x\in V}}
\rho u(x,t)m(x)=2/t\left\Vert u_{0}\right\Vert _{\infty}%
{\displaystyle\sum\limits_{x\in V}}
\rho u_{0}(x)m(x)\text{.}%
\]
Here we have used the mass conservation $\left\Vert \rho u(t)\right\Vert
_{1}=\left\Vert \rho u_{0}\right\Vert _{1},$ $t>0$. For the proof we refer
[33], [2]. Note that under the condition $(H_{1})$ the proof of mass
conservation is quite similar.

Lemma 2.3 is proved. $\square$

\textbf{Remark 2.1.} Let $u(x,t)$ be the solution of (1.1) and (1.2). Then
from Lemma 2.3 it follows that for $t>t_{0}$, where $t_{0}$ is large enough an
assumption (2.2) holds with $r=1$. Therefore, (2.17) holds as well, where for
any fixed $t>t_{0}(u_{0})$ we set%

\[
E_{2}(t):=%
{\displaystyle\sum\limits_{x\in V}}
\rho u^{2}(x,t)m(x)\text{, }D_{p}(t):=%
{\displaystyle\sum\limits_{x\thicksim y}}
\left\vert D_{y}u(x,t)\right\vert ^{p}\mathit{w(x,y)}\text{,}%
\]

\[
E_{1}(t)=%
{\displaystyle\sum\limits_{x\in V}}
\rho u(x,t)m(x)\text{.}%
\]

\textbf{Lemma 2.4}. \textit{Under the assumptions }$(H_{1})$\textit{ we have
for }$R>0$%

\[
\mathit{\gamma}^{1/A_{2}}\mathit{\Phi}^{\left(  -1\right)  }\mathit{(R)\leq
\Phi}^{(-1)}\mathit{(\gamma R)\leq\gamma}^{1/A_{1}}\mathit{\Phi}%
^{(-1)}\mathit{(R)}\text{, }\mathit{0<\gamma<1,}%
\]

\[
\mathit{\gamma}^{1/A_{2}}\mathit{\Phi}^{\left(  -1\right)  }\mathit{(R)\leq\Phi}^{(-1)}\mathit{(\gamma
R)\leq\gamma}^{1/A_{1}}\mathit{\Phi}^{(-1)}(R)\text{, }\mathit{\gamma\geq1}\text{,}%
\]
\textit{where}%

\[
\mathit{A}_{1}\mathit{=}\frac{\lambda+N+\alpha_{1}-p\alpha_{2}-p}{\left(
\lambda+N-p\alpha_{2}\right)  (p-1)}\mathit{,}\text{ }\mathit{A}_{2}%
\mathit{=}\frac{\lambda+N+\alpha_{2}-p\alpha_{1}-p}{\left(  \lambda
+N-p\alpha_{1}\right)  (p-1)}\mathit{>0}.%
\]
\textit{Proof}. We have for any $R>0$%

\[
\gamma^{p-\alpha_{1}}\omega(R)\leq\omega(\gamma R)\leq\gamma^{p-\alpha_{2}%
}\omega(R)\text{, }0<\gamma<1,
\]

\[
\gamma^{p-\alpha_{2}}\omega(R)\leq\omega(\gamma R)\leq\gamma^{p-\alpha_{1}%
}\omega(R),\text{ }\gamma\geq1\text{.}%
\]
Thus,%

\[
\gamma^{\left(  \lambda+N-p\alpha_{1}\right)  /p}\Psi(R)\leq\Psi(\gamma
R)\leq\gamma^{\left(  \lambda+N-p\alpha_{2}\right)  /p}\Psi(R),\text{
}0<\gamma<1\text{,}%
\]

\[
\gamma^{\left(  \lambda+N-p\alpha_{2}\right)  /p}\Psi(R)\leq\Psi(\gamma
R)\leq\gamma^{\left(  \lambda+N-p\alpha_{1}\right)  /p}\Psi(R),\text{ }%
\gamma\geq1,
\]

\[
\gamma^{p/\left(  \lambda+N-p\alpha_{2}\right)  }\Psi^{(-1)}(R)\leq\Psi
^{(-1)}(\gamma R)\leq\gamma^{p/\left(  \lambda+N-p\alpha_{1}\right)  }%
\Psi^{(-1)}(R)\text{, }0<\gamma<1,
\]

\[
\gamma^{p/\left(  \lambda+N-p\alpha_{1}\right)  }\Psi^{(-1)}(R)\leq\Psi
^{(-1)}(\gamma R)\leq\gamma^{p/\left(  \lambda+N-p\alpha_{2}\right)  }%
\Psi^{(-1)}(R)\text{, }\gamma\geq1.
\]
In the same way we have%

\[
\gamma^{A_{1}}\Phi(R)\leq\Phi(\gamma R)\leq\gamma^{A_{2}}\Phi(R)\text{,
}0<\gamma<1,
\]

\[
\gamma^{A_{2}}\Phi(R)\leq\Phi(\gamma R)\leq\gamma^{A_{1}}\Phi(R)\text{,
}\gamma\geq1.
\]
Finally,%

\[
\gamma^{1/A_{2}}\Phi^{\left(  -1\right)  }(R)\leq\Phi^{(-1)}(\gamma
R)\leq\gamma^{1/A_{1}}\Phi^{(-1)}(R)\text{, }0<\gamma<1\text{,}%
\]

\[
\gamma^{1/A_{1}}\Phi^{(-1)}(R)\leq\Phi^{(-1)}(\gamma R)\leq\gamma^{1/A_{2}%
}\Phi^{\left(  -1\right)  }(R),\text{ }\gamma\geq1.
\]

Lemma 2.4 is proved. $\square$

In the proofs of main results we will use the following classical recursive lemmas.

\textbf{Lemma 2.5.} (see   0 of [13]). \textit{Let }$\left\{
Y_{n}\right\}  _{n=1}^{\infty}$\textit{ be a sequence of equi-bounded positive
numbers satisfying recursive inequalities}%

\[
\mathit{Y}_{n}\mathit{\leq Cb}^{n}\mathit{Y}_{n+1}^{1-\alpha}\text{,}%
\]

\textit{where }$C,$\textit{ }$b>1$\textit{ and }$\alpha\in(0,1)$\textit{ are
given constants. Then}%

\[
\mathit{Y}_{0}\mathit{\leq(2Cb}^{(1-\alpha)/\alpha}\mathit{)}^{1/\alpha
}\text{.}%
\]

\textbf{Lemma 2.6}. (see [29] Chapter 2, Section 5) \textit{Let }$\left\{
Y_{n}\right\}  _{n=1}^{\infty}$\textit{ be a sequence of equi-bounded positive
numbers satisfying recursive inequalities}%

\[
\mathit{Y}_{n}\mathit{\leq Cb}^{n}\mathit{Y}_{n+1}^{1+\alpha}\text{,}%
\]
\textit{where }$C,$\textit{ }$b>1$\textit{ and }$\alpha>0$\textit{, are given
constants. Then, if }%

\[
\mathit{Y}_{0}\mathit{\leq C}^{-1/\alpha}\mathit{b}^{-1/\alpha^{2}}%
\]
\textit{then }$Y_{n}\rightarrow0$\textit{ when }$n\rightarrow\infty$\textit{.}

\begin{center}
\textbf{3 Proof of Theorem 1.1.}
\end{center}

We start with stating the next technical lemma, which was proven in [2] :

\textbf{Lemma 3.1.} \textit{Let \ }$q>0$\textit{, }$p>2$\textit{, }$h\geq0,u,v:$%
\textit{\ }$V\rightarrow \mathbb{R}$\textit{. Then for all }$x,y\in V$%

\[
\mathit{(}\left\vert D_{y}u(x)\right\vert ^{p-2}\mathit{D}_{y}\mathit{u(x)-}%
\left\vert D_{y}v(x)\right\vert ^{p-2}\mathit{D}_{y}\mathit{v(x))}
\]

\begin{equation}
\mathit{\times D}_{y}\mathit{(u(x)-v(x)-h)}_{+}^q\mathit{\geq C}\left\vert
D_{y}(u(x)-v(x)-h)_{+}^{\frac{q-1+p}{p}}\right\vert ^{p}\mathit{,} \tag{3.1}
\end{equation}
where $C = C(p,q) > 0$ depends only on $p$ and $q$.

\textit{Proof.} For any $\tau_{1}>\tau_{2}>0$ we define the standard nonnegative cut-off function $\eta(t)\in C^{1}(\mathbb{R})$ such that%

\[
\eta(t)=0,\text{ }t\geq\tau_{1};\text{ }\eta(t)=0\text{, }t\leq\tau_{2}\text{;
}0\leq\eta^{\prime}(t)\leq\frac{2}{\tau_{1}-\tau_{2}}\text{, }t\in
\mathbb{R}\text{.}%
\]
Multiplying the both sides of (1.1) by $g(x,t):=\eta(t)^{p}(u(x,t)-h)_{+}^q$,
and using integration by parts formula:%

\[%
{\displaystyle\sum\limits_{x,y\in V}}
\left\vert D_{y}u(x,t))\right\vert ^{p-2}\mathit{D}_{y}\mathit{u(x,t)g(x,t)}%
w(x,y)=
\]%
\[
=-\frac{1}{2}%
{\displaystyle\sum\limits_{x,y\in V}}
\left\vert D_{y}u(x,t))\right\vert ^{p-2}\mathit{D}_{y}\mathit{u(x,t)}%
D_{y}g(x,t)w(x,y),
\]
we obtain%

\[
\frac{1}{1+q}%
{\displaystyle\sum\limits_{x\in V}}
\rho(x)(u(x,t)-h)_{+}^{1+q}m(x)\eta(t)^{p}%
\]

\[
+\frac{1}{2}%
{\displaystyle\int\limits_{0}^{t}}
{\displaystyle\sum\limits_{x,y\in V}}
\left\vert D_{y}u(x,t)\right\vert ^{p-2}D_{y}u(x,\tau)D_{y}[(u(x,\tau
)-h)_{+}^q]w(x,y)\eta(t)^{p}d\tau
\]

\begin{equation}
=\frac{p}{1+q}%
{\displaystyle\int\limits_{0}^{t}}
{\displaystyle\sum\limits_{x\in V}}
\rho(x)(u(x,t)-h)_{+}^{1+q}\eta^{\prime}(\tau)\eta(\tau)^{p-1}m(x)d\tau.
\tag{3.2}%
\end{equation}
By the Lemma 3.1 with $v=0$ we have%

\[
\left\vert D_{y}u(x,t)\right\vert ^{p-2}D_{y}u(x,\tau)D_{y}[(u(x,\tau
)-h)_{+}^q]
\]

\begin{equation}
\geq C\left\vert D_{y}(u(x,\tau)-h)_{+}^\frac{q-1+p}{p}\right\vert ^{p}. \tag{3.3}%
\end{equation}
Thus, combining (3.2) and (3.3), we get%

\[
\frac{1}{1+q}%
{\displaystyle\sum\limits_{x\in V}}
\rho(x)(u(x,t)-h)_{+}^{1+q}m(x)\eta(t)^{p}%
\]

\[
+C%
{\displaystyle\int\limits_{0}^{t}}
{\displaystyle\sum\limits_{x,y\in V}}
\left\vert D_{y}(u(x,\tau)-h)_{+}^\frac{q-1+p}{p}\right\vert ^{p}\eta(\tau)^{p}w(x,y)d\tau
\]

\begin{equation}
\leq\frac{p}{1+q}%
{\displaystyle\int\limits_{0}^{t}}
{\displaystyle\sum\limits_{x\in V}}
\rho(x)(u(x,t)-h)_{+}^{1+q}\eta^{\prime}(\tau)\eta(\tau)^{p-1}m(x)d\tau.
\tag{3.4}%
\end{equation}
Using (3.4) and the properties of the cut-off function $\eta$, we have%

\[
\sup_{\tau_{1}<\tau<t}%
{\displaystyle\sum\limits_{x\in V}}
\rho(x)(u(x,\tau)-h)_{+}^{1+q}m(x)
\]

\[
+%
{\displaystyle\int\limits_{0}^{t}}
{\displaystyle\sum\limits_{x,y\in V}}
\left\vert D_{y}(u(x,\tau)-h)_{+}^\frac{q-1+p}{p}\right\vert ^{p}w(x,y)d\tau
\]

\begin{equation}
\leq C\frac{1}{\tau_{1}-\tau_{2}}%
{\displaystyle\int\limits_{0}^{t}}
{\displaystyle\sum\limits_{x\in V}}
\rho(x)(u(x,\tau)-h)_{+}^{1+q}m(x)d\tau.\tag{3.5}%
\end{equation}
Define :%

\[
0<\sigma_{1}<\sigma_{2}<1/2,t_{i}=t/2(1-\sigma_{2})+t2^{-i-1}(\sigma
_{2}-\sigma_{1}),
\]

\[
k_{i}=k(1-\sigma_{2})+k2^{-i}(\sigma_{2}-\sigma_{1})\tau_{1}=t_{i},\tau
_{2}=t_{i+1}%
\]
and $i=0,1,..,$ $f_{i}=(u-k_{i})_{+}^\frac{q-1+p}{p}.$ Then (3.5) will take the form%

\[
\sup_{t_{i}<\tau<t}%
{\displaystyle\sum\limits_{x\in V}}
\rho(x)f_{i}^{1+q}(x,\tau)m(x)+%
{\displaystyle\int\limits_{t_{i}}^{t}}
{\displaystyle\sum\limits_{x,y\in V}}
\left\vert D_{y}f_{i}^\frac{q-1+p}{p}(x,\tau)\right\vert ^{p}w(x,y)d\tau
\]

\begin{equation}
\leq C\frac{2^{i}}{t(\sigma_{2}-\sigma_{1})}%
{\displaystyle\int\limits_{t_{i+1}}^{t}}
{\displaystyle\sum\limits_{x\in V}}
\rho(x)f_{i+1}^{1+q}(x,\tau)m(x)d\tau. \tag{3.6}%
\end{equation}
By Lemma 2.3 with (see Remark 2.1 and (2.15)) $f=f_{i+1}$, $q=1$, $r=1$ we get%

\begin{equation}
\mathit{KE}_{q}(f_{i+1})\mathit{\leq\varepsilon D}(f_{i+1})\mathit{+C}%
_{\varepsilon}\mathit{E}_{1}^{p}(f_{i+1})\mathit{\Phi}_{1}^{(-1)}\left(
\frac{K}{E_{1}^{p-2}(f_{i+1})}\right)  \text{,} \tag{3.7}%
\end{equation}
where%

\[
\mathit{E}_{q}(f_{i+1}):=%
{\displaystyle\sum\limits_{x\in V}}
\rho(x)f_{i+1}^{2}(x,\tau)m(x),\text{ }%
\]%
\[
D(f_{i+1})=%
{\displaystyle\sum\limits_{x,y\in V}}
\left\vert D_{y}f_{i+1}(x,\tau)\right\vert ^{p}w(x,y).
\]
Thus, choosing in (3.7) $K=C2^{i}t^{-1}(\sigma_{2}-\sigma_{1})^{-1}$ and
integrating in time, we obtain%

\[
C\frac{2^{i}}{t(\sigma_{2}-\sigma_{1})}%
{\displaystyle\int\limits_{t_{i+1}}^{t}}
\mathit{E}_{2}(f_{i+1}(\tau))d\tau
\]

\[
\leq\mathit{\varepsilon%
{\displaystyle\int\limits_{t_{i+1}}^{t}}
D}(f_{i+1}(\tau))d\tau\mathit{+C}_{\varepsilon}\frac{B^{i}}{(\sigma_{2}%
-\sigma_{1})^{c}}%
\]

\begin{equation}
\times\sup_{t_{\infty}<\tau<t}\mathit{E}_{r}(f_{\infty}(\tau))\left(  \Phi
_{1}^{(-1)}\left(  \frac{1}{\left(  t-t_{\infty}\right)  \sup\limits_{t_{\infty}%
<\tau<t}E_{1}^{p-2}(f_{\infty}(\tau))}\right)  \right)  \text{.} \tag{3.8}%
\end{equation}
Here we have used Lemma 2.5 : $\mathit{\Phi}_{1}^{(-1)}(\gamma s)\leq
\gamma^{c}\mathit{\Phi}_{1}^{(-1)}(s)$, $c>1$, with $\gamma=C2^{i}(\sigma
_{2}-\sigma_{1})^{-1}$, $c=c(\alpha_{1},\alpha_{2},N,p)$. Now from (3.7) we have%

\[
Y_{i}:=\sup_{t_{i}<\tau<t}%
{\displaystyle\sum\limits_{x\in V}}
\rho(x)f_{i}^{2}(x,\tau)m(x)+%
{\displaystyle\int\limits_{t_{i}}^{t}}
{\displaystyle\sum\limits_{x,y\in V}}
\left\vert D_{y}f_{i}(x,\tau)\right\vert ^{p}w(x,y)d\tau
\]

\[
\leq\varepsilon Y_{i+1}+\mathit{C}_{\varepsilon}C\frac{B^{i}}{(\sigma
_{2}-\sigma_{1})^{c}}%
\]

\[
\times\sup_{t_{\infty}<\tau<t}\mathit{E}_{1}(f_{\infty}(\tau))\mathit{\Phi
}_{1}^{(-1)}\left(  \frac{1}{\left(  t-t_{\infty}\right)  \sup\limits_{t_{\infty
}<\tau<t}E_{1}^{p-2}(f_{\infty}(\tau))}\right)  .
\]
Choosing $\varepsilon B<1$ by the standard iteration (see also Lemma 2.6)), we get%

\[
\sup_{t_{0}<\tau<t}%
{\displaystyle\sum\limits_{x\in V}}
\rho(x)f_{0}^{2}(x,\tau)m(x)\leq\frac{C}{(\sigma_{2}-\sigma_{1})^{c}}%
\]

\begin{equation}
\times\sup_{t_{\infty}<\tau<t}\mathit{E}_{1}(f_{\infty}(\tau))\left(  \Phi
_{1}^{(-1)}\left(  \frac{1}{\left(  t-t_{\infty}\right)  \sup\limits_{t_{\infty}%
<\tau<t}E_{1}^{p-2}(f_{\infty}(\tau))}\right)  \right)  . \tag{3.9}%
\end{equation}
Set%

\[
\sigma_{1}=\delta2^{-n-1},\sigma_{2}=\delta2^{-n},t_{0}\rightarrow
t/2(1-\delta2^{-n-1}),
\]

\[
t_{\infty}\rightarrow t/2(1-\delta2^{-n}),k_{0}\rightarrow k(1-\delta
2^{-n-1}),\overline{k}_{n}=(k_{n}+k_{n+1})/2.
\]
Then, applying (3.9) with $t_{0}=t_{n+1},t_{\infty}=t_{n},k_{0}=\overline
{k}_{n},k_{\infty}=k_{n}$, we obtain%

\[
I_{n}:=\sup_{t_{n}<\tau<t}%
{\displaystyle\sum\limits_{x\in V}}
\rho(x)\left(  u(x,\tau)-\overline{k}_{n}\right)  _{+}m(x)
\]

\begin{equation}
\leq C2^{cn}k^{-1}I_{n-1}^{2}\left(  \Phi_{1}^{(-1)}\left(  \frac{1}%
{tI_{n-1}^{p-2}}\right)  \right)  . \tag{3.10}%
\end{equation}
By iterative Lemma 2.7 we conclude that $I_{n}\rightarrow0$ as $n\rightarrow
\infty$, and therefore $\left\Vert u(t)\right\Vert _{\infty}\leq k$, provided%

\begin{equation}
k^{-1}I_{0}\left(  \Phi_{1}^{(-1)}\left(  \frac{1}{tI_{0}^{p-2}}\right)
\right)  \leq\delta\text{,} \tag{3.11}%
\end{equation}
where $\delta$ is a sufficiently small number which depends only on parameters
of the problem. As a matter of fact that the function $s\Phi_{1}^{(-1)}(s)$
is increasing, we conclude that (3.11) is in force if the following inequality holds%

\begin{equation}
k^{-1}M(t)\left(  \Phi_{1}^{(-1)}\left(  \frac{1}{tM^{p-2}(t)}\right)
\right)  \leq\delta\text{,} \tag{3.12}%
\end{equation}
where%

\[
M(t)=\sup_{0\leq\tau\leq t}%
{\displaystyle\sum\limits_{x\in V}}
\rho(x)u(x,\tau)m(x)\text{.}%
\]
Next, we need the mass conservation: for any $t>0$%

\begin{equation}
M(t)=M(0). \tag{3.13}%
\end{equation}
Choose now $k$ as follows%

\begin{equation}
k^{-1}M(0)\left(  \Phi_{1}^{(-1)}\left(  \frac{1}{tM^{p-2}(0)}\right)
\right)  =\delta. \tag{3.14}%
\end{equation}
Finally we transform (3.14) to the form (1.7)$.$ Namely, by the definition of
the functions $\Phi,\Phi_{1}$, $\Psi$ we consistently have%

\[
\left(  \Phi_{1}^{(-1)}\left(  \frac{1}{tM^{p-2}(0)}\right)  \right)
=\frac{\delta k}{M(0)}\iff
\]

\[
\frac{1}{tM^{p-2}(0)}=\Phi_{1}\left(  \frac{\delta k}{M(0)}\right)
=\frac{M(0)}{\delta k}\Phi^{(-1)}\left(  \frac{\delta k}{M(0)}\right)  \iff
\]

\[
\frac{\delta k}{tM^{p-1}(0)}=\Phi^{(-1)}\left(  \frac{\delta k}{M(0)}\right)
\text{ }\iff\frac{\delta k}{M(0)}=\Phi\left(  \frac{\delta k}{tM^{p-1}%
(0)}\right)
\]

\[
=\left(  \left(  \frac{\delta k}{tM^{p-1}(0)}\right)  \left[  \omega\left(
\Psi^{(-1)}\left(  \frac{tM^{p-1}(0)}{\delta k}\right)  ^{1/p}\right)
\right]  \right)  ^{1/(p-1)}\iff
\]

\[
\iff\omega\left(  \Psi^{(-1)}\left(  \frac{tM^{p-1}(0)}{\delta k}\right)
^{1/p}\right)  =\delta^{p-2}k^{p-2}t
\]

\[
\iff\Psi^{(-1)}\left[  \left(  \frac{tM^{p-1}(0)}{\delta k}\right)
^{1/p}\right]  =\omega^{(-1)}(\delta^{(p-2)}k^{p-2}t)
\]

\[
\iff\left(  \frac{tM^{p-1}(0)}{\delta k}\right)  ^{1/p}=\Psi\left(
\omega^{(-1)}(\delta^{(p-2)}k^{p-2}t)\right)  \iff\text{ }\delta^{p-2}%
k^{p-2}t
\]

\[
\times\omega^{(-1)}(\delta^{p-2}k^{p-2}t)^{(N-p)(p-1)/p}%
\]

\[
\left(  \left(  \omega^{(-1)}(\delta^{p-2}k^{p-2}t\right)  ^{N-p}\left(
\delta k\right)  ^{p-1}tM(0)\right)  ^{(p-1)/p}=1\text{,}%
\]

\[
\iff\left(  \omega^{(-1)}(\delta^{p-2}k^{p-2}t\right)  ^{N-p}\left(  \delta
k\right)  ^{p-1}tM(0)=1.
\]
Since $\left\Vert u(t)\right\Vert _{\infty}\leq k$, by Lemma 2.5 we arrive at
the desired result%

\[
\left(  \omega^{(-1)}(\left\Vert u(t)\right\Vert _{\infty}^{p-2}t)\right)
^{N-p}\left\Vert u(t)\right\Vert _{\infty}^{p-1}tM(0)\leq C\text{.}%
\]
Theorem 1.1 is proved. $\square$

\begin{center}
\textbf{4 Proof of Theorem 1.2. }
\end{center}

\bigskip Proceeding as in the proof of Theorem 1.1 we use the Caccioppoli
inequality (see Lemma 3.1.): for any $\tau_{1}>\tau_{2}>0$ we define the
standard nonnegative cut-off function $\eta(t)\in C^{1}(\mathbb{R})$ such that%

\[
\eta(t)=0,\text{ }t\geq\tau_{1};\text{ }\eta(t)=0\text{, }t\leq\tau_{2}\text{;
}0\leq\eta^{\prime}(t)\leq\frac{2}{\tau_{1}-\tau_{2}}\text{, }t\in
\mathbb{R}\text{.}%
\]
Multiplying the both sides of (1.1) by $g(x,t):=\eta(t)^{p}(u(x,t)-h)_{+}%
^{\theta}$, $\theta>0$, and using the integration by parts formula:%

\[%
{\displaystyle\sum\limits_{x,y\in V}}
\left\vert D_{y}u(x,t))\right\vert ^{p-2}\mathit{D}_{y}\mathit{u(x,t)g(x,t)}%
w(x,y)=
\]%
\[
=-\frac{1}{2}%
{\displaystyle\sum\limits_{x,y\in V}}
\left\vert D_{y}u(x,t))\right\vert ^{p-2}\mathit{D}_{y}\mathit{u(x,t)D}%
_{y}u(x,t)D_{y}g(x,t)w(x,y),
\]
we obtain%

\[
\frac{1}{1+\theta}%
{\displaystyle\sum\limits_{x\in V}}
\rho(x)(u(x,t)-h)_{+}^{1+\theta}m(x)\eta(t)^{p}%
\]

\[
+\frac{1}{2}%
{\displaystyle\int\limits_{0}^{t}}
{\displaystyle\sum\limits_{x,y\in V}}
\left\vert D_{y}u(x,t)\right\vert ^{p-2}D_{y}u(x,\tau)D_{y}[(u(x,\tau
)-h)_{+}^{\theta}]w(x,y)\eta(t)^{p}d\tau
\]

\begin{equation}
=\frac{p}{1+\theta}%
{\displaystyle\int\limits_{0}^{t}}
{\displaystyle\sum\limits_{x\in V}}
\rho(x)(u(x,t)-h)_{+}^{1+\theta}\eta^{\prime}(\tau)\eta(\tau)^{p-1}%
m(x)d\tau. \tag{4.1}%
\end{equation}
By the Lemma 3.1 with $v=0$ we have%

\begin{equation}
	\left\vert D_{y}u(x,t)\right\vert ^{p-2}D_{y}u(x,\tau)D_{y}[(u(x,\tau)-h)_{+}^{\theta}] \geq C_{0}\left\vert D_{y}(u(x,\tau)-h)_{+}^{(p+\theta-1)/p}\right\vert ^{p}.
	\tag{4.2}%
\end{equation}
Thus, combining (4.1) and (4.2), we get%

\[
\frac{1}{1+\theta}%
{\displaystyle\sum\limits_{x\in V}}
\rho(x)(u(x,\tau)-h)_{+}^{1+\theta}m(x)\eta(t)^{p}%
\]

\[
+C_{0}%
{\displaystyle\int\limits_{0}^{t}}
{\displaystyle\sum\limits_{x,y\in V}}
\left\vert D_{y}(u(x,\tau)-h)_{+}^{(p+\theta-1)/p}\right\vert ^{p}\eta
(\tau)^{p}w(x,y)d\tau
\]

\begin{equation}
\leq\frac{p}{1+\theta}%
{\displaystyle\int\limits_{0}^{t}}
{\displaystyle\sum\limits_{x\in V}}
\rho(x)(u(x,\tau)-h)_{+}^{1+\theta}\eta^{\prime}(\tau)\eta(\tau)^{p-1}%
m(x)d\tau. \tag{4.3}%
\end{equation}
Using (4.3) and the properties of the cut-off function $\eta$, we have%

\[
\sup_{\tau_{1}<\tau<t}%
{\displaystyle\sum\limits_{x\in V}}
\rho(x)(u(x,\tau)-h)_{+}^{1+\theta}m(x)
\]

\[
+%
{\displaystyle\int\limits_{0}^{t}}
{\displaystyle\sum\limits_{x,y\in V}}
\left\vert D_{y}(u(x,\tau)-h)_{+}^{(p+\theta-1)/p}\right\vert ^{p}w(x,y)d\tau
\]

\begin{equation}
\leq C\frac{1}{\tau_{1}-\tau_{2}}%
{\displaystyle\int\limits_{0}^{t}}
{\displaystyle\sum\limits_{x\in V}}
\rho(x)(u(x,\tau)-h)_{+}^{1+\theta}m(x)d\tau. \tag{4.4}%
\end{equation}
Define :%

\[
0<\sigma_{1}<\sigma_{2}<1/2,t_{i}=t/2(1-\sigma_{2})+t2^{-i-1}(\sigma
_{2}-\sigma_{1}),
\]

\[
k_{i}=k(1-\sigma_{2})+k2^{-i}(\sigma_{2}-\sigma_{1})\tau_{1}=t_{i},\tau
_{2}=t_{i+1}\text{,}%
\]

\[
i=0,1,..,f_{i}=(u-k_{i})_{+}^{s},s=(p+\theta-1)/p\text{.}%
\]
Let%

\[
0<\nu<\theta,a=(1+\nu)p/(p+\theta-1),b=(1+\theta)p/(p+\theta-1).
\]
Then (4.4) reads%

\[
\sup_{t_{i}<\tau<t}%
{\displaystyle\sum\limits_{x\in V}}
\rho(x)f_{i}^{b}(x,\tau)m(x)+%
{\displaystyle\int\limits_{t_{i}}^{t}}
{\displaystyle\sum\limits_{x,y\in V}}
\left\vert D_{y}f_{i}(x,\tau)\right\vert ^{p}w(x,y)d\tau
\]

\begin{equation}
\leq C\frac{2^{i}}{t(\sigma_{2}-\sigma_{1})}%
{\displaystyle\int\limits_{0}^{t}}
{\displaystyle\sum\limits_{x\in V}}
\rho(x)f_{i+1}^{r}(x,\tau)m(x)d\tau. \tag{4.5}%
\end{equation}
\bigskip Applying the H\"{o}lder inequality, and (1.3) with $f=f_{i+1}$, we
get \
\[%
{\displaystyle\sum\limits_{x\in V}}
\rho(x)f_{i+1}^{b}m(x)\leq\left(
{\displaystyle\sum\limits_{x\in V}}
f_{i+1}^{p^{\ast}}m(x)\right)  ^{(b-r)/(p^{\ast}-r)}%
\]

\[
\times\left(
{\displaystyle\sum\limits_{x\in V}}
\rho(x)^{(p^{\ast}-r)/(p^{\ast}-b)}f_{i+1}^{r}m(x)\right)  ^{(p^{\ast
}-b)/(p^{\ast}-r)}%
\]

\begin{equation}
\leq C\left(
{\displaystyle\sum\limits_{x,y\in V}}
\left\vert D_{y}f_{i+1}(x,\tau)\right\vert ^{p}w(x,y)\right)  ^{A}\left(
{\displaystyle\sum\limits_{x\in V}}
\rho(x)f_{i+1}^{r}m(x)\right)  ^{B}, \tag{4.6}%
\end{equation}
where%

\[
A=\frac{b-r}{p^{\ast}-r}\frac{p^{\ast}}{p}=\frac{N(\theta-\nu)}{\lambda
+N(\theta-\nu)+p\nu}<1,
\]

\[
B=\frac{\lambda+p\theta}{\lambda+N(\theta-\nu)+p\nu}.%
\]

Note that in (4.6) we have used the elementary fact: $\rho(x)\leq1$ and
$(p^{\ast}-r)/(p^{\ast}-b)>1$ imply $\rho(x)^{(p^{\ast}-r)/(p^{\ast}-b)}%
\leq\rho(x)$. Next, integrating (4.6) in time and using H\"{o}lder's
inequality, we obtain%

\[%
{\displaystyle\int\limits_{t_{i+1}}^{t}}
d\tau%
{\displaystyle\sum\limits_{x\in V}}
\rho(x)f_{i+1}^{b}m(x)
\]

\[
\leq C%
{\displaystyle\int\limits_{t_{i+1}}^{t}}
d\tau\left(
{\displaystyle\sum\limits_{x,y\in V}}
\left\vert D_{y}f_{i+1}(x,\tau)\right\vert ^{p}w(x,y)\right)  ^{A}%
\]

\[
\times\left(
{\displaystyle\sum\limits_{x\in V}}
\rho(x)f_{i+1}^{r}m(x)\right)  ^{B}%
\]

\[
\leq C\left[
{\displaystyle\int\limits_{t_{i+1}}^{t}}
d\tau\left(
{\displaystyle\sum\limits_{x,y\in V}}
\left\vert D_{y}f_{i+1}(x,\tau)\right\vert ^{p}w(x,y)\right)  \right]  ^{A}%
\]

\[
\times t^{1-A}\left(  \sup_{t_{i+1}<\tau<t}%
{\displaystyle\sum\limits_{x\in V}}
\rho(x)f_{i+1}^{r}(x,\tau)m(x)\right)  ^{B}.
\]
Thus, applying the Young inequality with $\varepsilon>0$, we get%

\[
C\frac{2^{i}}{(\sigma_{2}-\sigma_{1})t}%
{\displaystyle\int\limits_{t_{i+1}}^{t}}
d\tau%
{\displaystyle\sum\limits_{x\in V}}
\rho(x)f_{i+1}^{b}m(x)
\]

\[
\leq\varepsilon^{\frac{1}{A}}%
{\displaystyle\int\limits_{t_{i+1}}^{t}}
d\tau%
{\displaystyle\sum\limits_{x,y\in V}}
\left\vert D_{y}f_{i+1}(x,\tau)\right\vert ^{p}w(x,y)
\]

\begin{equation}
+\left(  \frac{C2^{i}}{\varepsilon(\sigma_{2}-\sigma_{1})t^{A}}\right)
^{1/(1-A)}\times\left(  \sup_{t_{i+1}<\tau<t}%
{\displaystyle\sum\limits_{x\in V}}
\rho(x)f_{i+1}^{r}m(x)\right)  ^{B/(1-A)}\text{.}\tag{4.7}%
\end{equation}
Hence (4.5)-(4.7) yield:%

\[
Y_{i}:=\sup_{t_{i}<\tau<t}%
{\displaystyle\sum\limits_{x\in V}}
\rho(x)\left(  u(x,\tau)-k_{i}\right)  _{+}^{1+\theta}m(x)
\]

\[
+%
{\displaystyle\int\limits_{t_{i}}^{t}}
{\displaystyle\sum\limits_{x,y\in V}}
\left\vert D_{y}f_{i}(x,\tau)\right\vert ^{p}w(x,y)d\tau
\]

\[
\leq\varepsilon^{1/A}Y_{i+1}+\left(  \frac{C2^{i}}{\varepsilon(\sigma
_{2}-\sigma_{1})t^{A}}\right)  ^{1/(1-A)}%
\]

\begin{equation}
\times\left(  \sup_{t_{i+1}<\tau<t}%
{\displaystyle\sum\limits_{x\in V}}
\rho(x)f_{i+1}^{r}m(x)\right)  ^{B/(1-A)}. \tag{4.8}%
\end{equation}
Note here that $Y_{i}\leq%
{\displaystyle\sum\limits_{x\in V}}
\rho(x)u_{0}(x)^{1+\theta}m(x)$ for all $i$. Iterating (4.8), we get by induction for all $i$%

\[
Y_{0}\leq\varepsilon^{i/(1-A)}Y_{i}+\left(  \frac{C}{\varepsilon(\sigma
_{2}-\sigma_{1})t^{A}}\right)  ^{1/(1-A)}\text{.}%
\]%
\[
\times\left(  \sup_{t_{i+1}<\tau<t}%
{\displaystyle\sum\limits_{x\in V}}
\rho(x)f_{i+1}^{r}m(x)\right)  ^{B/(1-A)}%
{\displaystyle\sum\limits_{j=0}^{i}}
(2\varepsilon)^{j/(1-A)}.%
\]
Next, choosing $\varepsilon:2\varepsilon<1$ and letting $i\rightarrow\infty$, we arrive at%

\[
\sup_{t_{0}<\tau<t}%
{\displaystyle\sum\limits_{x\in V}}
\rho(x)\left(  u(x,\tau)-k_{0}\right)  _{+}^{1+\theta}m(x)\leq Y_{0}%
\leq\frac{C}{(\sigma_{2}-\sigma_{1})^{1/(1-A)}t^{A/(1-A)}}%
\]

\begin{equation}
\leq\left(  \sup_{t_{\infty}<\tau<t}%
{\displaystyle\sum\limits_{x\in V}}
\rho(x)\left(  u(x,\tau)-k_{\infty}\right)  _{+}^{1+\nu}m(x)\right)
^{B/(1-A)}. \tag{4.9}%
\end{equation}
Set%

\[
\sigma_{1}=\delta2^{-n-1},\sigma_{2}=\delta2^{-n},t_{0}\rightarrow
t/2(1-\delta2^{-n-1}),
\]

\[
t_{\infty}\rightarrow t/2(1-\delta2^{-n}),k_{0}\rightarrow k(1-\delta
2^{-n-1}),\overline{k}_{n}=(k_{n}+k_{n+1})/2\text{.}%
\]
Applying the last inequality with $t_{0}=t_{n+1},t_{\infty}=t_{n}%
,k_{0}=\overline{k}_{n},k_{\infty}=k_{n}$, and taking into account that%

\[
\frac{A}{1-A}=\frac{N(\theta-\nu)}{\lambda+p\nu},\text{ }\frac{B}{1-A}%
=1+\frac{p(\theta-\nu)}{\lambda+p\nu},
\]
we conclude%

\[
I_{n}:=\sup_{t_{n}<\tau<t}%
{\displaystyle\sum\limits_{x\in V}}
\rho(x)\left(  u(x,\tau)-\overline{k}_{n}\right)  _{+}^{1+\nu}m(x)
\]

\[
\leq2^{cn}k^{-(\theta-\nu)}t^{-N(\theta-\nu)/(\lambda+p\nu)}I_{n-1}%
^{1+p(\theta-\nu)/(\lambda+p\nu)}.
\]
Iterating this inequality (see Lemma 2.5.), we deduce that $I_{n}\rightarrow0$
as $n\rightarrow\infty$ and therefore $\left\Vert u(t)\right\Vert _{\infty
}\leq2k$ provided%

\begin{equation}
k^{-(\theta-\nu)}t^{-N(\theta-\nu)/(\lambda+p\nu)}I_{0}^{p(\theta
-\nu)/(\lambda+p\nu)}\leq\epsilon. \tag{4.10}%
\end{equation}
Here $\epsilon$ is a sufficiently small enough number depending on the $N,p$ .
Then choosing%

\[
k=\epsilon^{-1/(\theta-\nu)}t^{-N/(\lambda+p\nu)}\left(  \sup_{t_{/4}<\tau<t}%
{\displaystyle\sum\limits_{x\in V}}
\rho(x)u(x,\tau)^{1+\nu}m(x)\right)  ^{p/(\lambda+p\nu)},
\]
we derive%

\begin{equation}
\left\Vert u(t)\right\Vert _{\infty}\leq Ct^{-N/(\lambda+p\nu)}\left(
\sup_{t_{/4}<\tau<t}%
{\displaystyle\sum\limits_{x\in V}}
\rho(x)u(x,\tau)^{1+\nu}m(x)\right)  ^{p/(\lambda+p\nu)}.\tag{4.11}%
\end{equation}
Our final task is to prove the integral form of universal bound:%

\begin{equation}%
{\displaystyle\sum\limits_{x\in V}}
\rho(x)u(x,t)^{1+\nu}m(x)\leq Ct^{-(1+\nu)/(p-2)}. \tag{4.12}%
\end{equation}
Then combining (4.11) and (4.12), we obtain%

\[
\left\Vert u(t)\right\Vert _{\infty}\leq Ct^{-N/(\lambda+p\nu)}\sup
_{t_{/4}<\tau<t}\tau^{-(1+\nu)p/(p-2)(\lambda+p\nu)}\leq Ct^{-1/(p-2)},
\]
as desired. Let us prove now (4.12). Multiplying both sides of the equation
(1.1) by $u^{\nu}(x,t)$ and integrating by parts, we obtain%

\begin{equation}
\frac{d}{dt}%
{\displaystyle\sum\limits_{x\in V}}
\rho(x)u(x,t)^{1+\nu}m(x)\leq-C%
{\displaystyle\sum\limits_{x,y\in V}}
\left\vert D_{y}u^{(p+\nu-1)/p}(x,\tau)\right\vert ^{p}w(x,y). \tag{4.13}%
\end{equation}
Let $\alpha=(p+\nu-1)/p$ and $v:=u^{\alpha}$. Applying H\"{o}lder's and
Sobolev's inequalities, we have%

\[%
{\displaystyle\sum\limits_{x\in V}}
\rho(x)v^{\beta}m(x)\leq\left(
{\displaystyle\sum\limits_{x\in V}}
\rho(x)v^{p^{\ast}}m(x)\right)  ^{\beta/p^{\ast}}%
\]

\[
\times\left(
{\displaystyle\sum\limits_{x\in V}}
\rho(x)^{p^{\ast}/(p^{\ast}-\beta)}m(x)\right)  ^{(p^{\ast}-\beta
)/p^{\ast}}%
\]

\begin{equation}
\leq C\left(
{\displaystyle\sum\limits_{x,y\in V}}
\left\vert D_{y}v(x,\tau)\right\vert ^{p}w(x,y)\right)  ^{\beta/p}\left(
{\displaystyle\sum\limits_{x\in V}}
\rho(x)^{p^{\ast}/(p^{\ast}-\beta)}m(x)\right)  ^{(p^{\ast}-\beta
)/p^{\ast}}. \tag{4.14}%
\end{equation}
Note that by assumption (1.4)%

\begin{equation}%
{\displaystyle\sum\limits_{x\in V}}
\rho(x)^{p^{\ast}/(p^{\ast}-\beta)}m(x)=%
{\displaystyle\sum\limits_{x\in V}}
\rho(x)^{N(p-1+\nu)/(\lambda+p\nu)}m(x)<\infty\tag{4.15}%
\end{equation}
for $\nu$ large enough. Therefore, denoting%

\[
E(t):=%
{\displaystyle\sum\limits_{x\in V}}
\rho(x)u(x,t)^{1+\nu}m(x)
\]
and combining (4.14), (4.15), we obtain%

\[
\frac{d}{dt}E(t)\leq-CE(t)^{(p+\nu-1)/(1+\nu)}.
\]
Integrating this inequality, we arrive at the desired estimate (4.10).
$\square$

\bigskip\textbf{Acknowledgements}

The author would like to express sincere gratitude to the anonymous reviewers
for their insightful comments and constructive suggestions, which
significantly contributed to improving the clarity and rigor of this
manuscript. The author also acknowledges the support and encouragement
provided by colleagues and mentors throughout the research process.

\bigskip\textbf{Data availability statement}

No datasets were generated or analysed during the current study.

\bigskip\textbf{Disclosure statement}

No potential conflict of interest was reported by the author.

\bigskip\bigskip

\section*{References{\protect\small \renewcommand{\labelenumi}{\theenumi.} }}

\begin{enumerate}
	\item {\small Andres, S., Deuschel, J.D., and Slowik, M., Heat kernel estimates for random walks with degenerate weights. Electron. J. Probab. \textbf{21}(33) (2016), 1–21. DOI: 10.1214/16-EJP4382.}
	
	\item {\small Andreucci, D. and Tedeev, A.F., Asymptotic estimates for the p-Laplacian on infinite graphs with decaying initial data. Potential Anal. \textbf{53} (2020), 677–699. DOI: 10.1007/s11118-019-09784-w.}
	
	\item {\small Andreucci, D. and Tedeev, A.F., Asymptotic Properties of Solutions to the Cauchy Problem for Degenerate Parabolic Equations with Inhomogeneous Density on Manifolds. Milan J. Math. \textbf{89} (2021), 295–327. DOI: 10.1007/s00032-021-00335-w.}
	
	\item {\small Andreucci, D. and Tedeev, A.F., Optimal decay rate for degenerate parabolic equations on noncompact manifolds. Methods Appl. Anal. \textbf{22}(4) (2015), 359–376. DOI: 10.4310/MAA.2015.V22.N4.A2.}
	
	\item {\small Andreucci, D. and Tedeev, A.F., Universal bounds at the blow-up time for nonlinear parabolic equations. Adv. Differ. Equ. \textbf{10}(1) (2005), 89–120. DOI: 10.57262/ade/1355867897.}
	
	\item {\small Biagi, S., Meglioli, G., and Punzo, F., A Liouville theorem for elliptic equations with a potential on infinite graphs. Calc. Var. \textbf{63}, 165 (2024). DOI: 10.1007/s00526-024-02768-8.}
	
	\item {\small Bianchi, D., G. Setti, and R. Wojciechowski, The generalized porous medium equation on graphs: existence and uniqueness of solutions with $l^1$ data. Preprint (2022). arXiv: 2112.01733.}
	
	\item {\small Cardoso, D.M. and Pinheiro, S.J., Spectral Bounds for the k-Regular Induced Subgraph Problem. In: Bebiano, N. (ed.) Appl. Comput. Matrix Anal. MAT-TRIAD 2015. Springer Proc. Math. Stat. \textbf{192} (2017). DOI: 10.1007/978-3-319-49984-0\_7.}
	
	\item {\small Chung, F., Grigor'yan, A., and Yau, S.T., Higher Eigenvalues and Isoperimetric Inequalities on Riemannian manifolds and graphs. Comm. Anal. Geom. \textbf{8}(5) (2000), 969–1026. DOI: 10.1007/s000390050070.}
	
	\item {\small Chung, S.-Y. and Choi, M.-J., Blow-up Solutions and Global Solutions to Discrete p-Laplacian Parabolic Equations. Abstr. Appl. Anal. (2014), Art. ID 351675, 11 pp. DOI: 10.1155/2014/351675.}
	
	\item {\small Coulhon, T. and Grigor'yan, A., Random Walks on Graphs with Regular Volume Growth. Geom. Funct. Anal. \textbf{8} (1998), 656–701. DOI: 10.1007/s000390050070.}
	
	\item {\small D. Bianchi, A. G. Setti, and R. K. Wojciechowski, The generalized porous medium equation on graphs: existence and uniqueness of solutions with $l^1$ data. Preprint (2022). arXiv: 2112.01733.}
	
	\item {\small DiBenedetto, E., Degenerate Parabolic Equations. Springer-Verlag, New York (1993).}
	
	\item {\small Dzagoeva, L.F. and Tedeev, A.F., Asymptotic behavior of the solutions of doubly degenerate parabolic equations with inhomogeneous density. Vladikavkaz Math. J. \textbf{24}(3) (2022), 78–86. DOI: 10.46698/p6936-3163-2954-s.}
	
	\item {\small Elmoataz, A., Toutain, M., and Tenbrinck, D., On the p-Laplacian and $\infty$-Laplacian on graphs with applications in image and data processing. SIAM J. Imaging Sci. \textbf{8}(4) (2015), 2412–2451. DOI: 10.1137/15M1022793.}
	
	\item {\small Galaktionov, V.A., Kamin, S., Kersner, R., and Vazquez, J.L., Intermediate Asymptotics for Inhomogeneous Nonlinear Heat Conduction. J. Math. Sci. \textbf{120}(3) (2004), 1277–1294. \\ DOI: 10.1023/B:JOTH.0000016049.94192.aa.}
	
	\item {\small Grigor'yan, A., Analysis and Geometry on Graphs Part 1: Laplace operator on weighted graphs. Tsinghua University (2012).}
	
	\item {\small Grillo, G., Muratori, M., and Punzo, F., On the asymptotic behavior of solutions to the fractional porous medium equation with variable density. Discrete Contin. Dyn. Syst. \textbf{35}(12) (2015), 5927–5962. DOI: 10.3934/dcds.2015.35.5927.}
	
	\item {\small Guedda, M., Hihorst, D., and Peletier, M.A., Disappearing Interfaces in Nonlinear Diffusion. Adv. Math. Sci. Appl. \textbf{7} (1997), 695–710.}
	
	\item {\small Holopainen, I. and Soardi, P.M., p-harmonic functions on graphs and manifolds. Manuscripta Math. \textbf{94} (1997), 95–110. DOI: 10.1007/BF02677841.}
	
	\item {\small Holopainen, I. and Soardi, P., A strong Liouville theorem for p-harmonic functions on graphs. Ann. Acad. Sci. Fenn. Math. \textbf{22} (1996).}
	
	\item {\small Hua, B. and Mugnolo, D., Time regularity and long time behavior of parabolic p-Laplace equations on infinite graphs. J. Differ. Equ. \textbf{259} (2015), 6162–6190. DOI: 10.1016/j.jde.2015.07.018.}
	
	\item {\small Kamin, S. and Rosenau, P., Propagation of thermal waves in an inhomogeneous medium. Commun. Pure Appl. Math. \textbf{34} (1981), 831–852.}
	
	\item {\small Kamin, S. and Rosenau, P., Nonlinear diffusion in finite mass medium. Commun. Pure Appl. Math. \textbf{35} (1982), 113–127.}
	
	\item {\small Kamin, S., Reyes, G., and Vazquez, J.L., Long time behavior for the inhomogeneous PME in a medium with rapidly decaying density. Discrete Contin. Dyn. Syst. \textbf{26}(2) (2010), 521–549. DOI: 10.3934/dcds.2010.26.521.}
	
	\item {\small Keller, M. and Lenz, D., Unbounded Laplacians on graphs: basic spectral properties and the heat equation. Math. Model. Nat. Phenom. 5 (2010), 198--224. MR 2662456 Zbl 1207.47032}
	
	\item {\small Keller, M., Lenz, D., and Wojciechowski, R.K. (eds.), Analysis and Geometry on Graphs and Manifolds. Cambridge University Press (2020).}
	
	\item {\small Keller, M., Lenz, D., Schmidt, M., and Wirth, M., Diffusion determines the recurrent graph, Advances in Mathematics, Volume 269, (2015), Pages 364--398, ISSN 0001--8708, \\ https://doi.org/10.1016/j.aim.2014.10.003.}
	
	\item {\small Ladyzhenskaya, O.A., Solonnikov, V.A., and Ural'ceva, N.N., Linear and Quasi-Linear Equations of Parabolic Type. Transl. Math. Monogr. \textbf{23}, AMS (1968).}
	
	\item {\small Lee, Y.S. and Chung, S.Y., Extinction and positivity of solutions of the p-Laplacian evolution equation on networks. J. Math. Anal. Appl. \textbf{386} (2012), 581–592. DOI: 10.1016/j.jmaa.2011.08.023.}
	
	\item {\small Masamune, J., A Liouville property and its application to the Laplacian of an infinite graph.
		In M. Kotani, H. Naito, and T. Tate (eds.), Spectral analysis in geometry and number
		theory. Papers from the International Conference on the occasion of Toshikazu Sunada's
		60th birthday held at Nagoya University, Nagoya, August 6--10, (2007). Contemporary
		Mathematics 484. Amer. Math. Soc., Providence (RI), (2009), pp. 103--115. MR 1500141
		Zbl 1181.58018}
	
	\item {\small Mourrat, J.-C. and Otto, F., Anchored Nash inequalities and heat kernel bounds for static and dynamic degenerate environments. J. Funct. Anal. \textbf{270} (2016), 201–228. DOI: 10.1016/j.jfa.2015.09.020.}
	
	\item {\small Mugnolo, D., Parabolic theory of the discrete p-Laplace operator. Nonlinear Anal. \textbf{87} (2013), 33–60. DOI: 10.1016/j.na.2013.04.002.}
	
	\item {\small Ostrovskii, M.I., Sobolev spaces on graphs. Quest. Math. \textbf{28}(4) (2005), 501–523. DOI: 10.2989/16073600509486144.}
	
	\item {\small Pinchover, Y. and Tintarev, K. On positive solutions of p-Laplacian-type equations. Analysis, Partial Differential Equations and Applications – The Vladimir Maz’ya Anniversary Volume, Operator Theory: Advances and Applications, vol. 193, Birkhäuser, (2009), pp. 225–236 }
	
	\item {\small Pinchover, Yehuda. Large time behavior of the heat kernel. Journal of Functional Analysis. (2004), 206. 191-209. 10.1016/S0022-1236(03)00110-1. }
	
	\item {\small Reyes, G. and Vazquez, J.L., Long Time Behavior for the Inhomogeneous PME in a Medium with Slowly Decaying Density. Commun. Pure Appl. Anal. \textbf{8}(2) (2009), 493–508. DOI: 10.3934/cpaa.2009.8.493.}
	
	\item {\small Saloff-Coste, L. Aspects of Sobolev-Type Inequalities, Volume 289 of London.Mathematical Society Lecture Note Series. Cambridge University Press, Cambridge (2002)}
	
	\item {\small Saloff-Coste, L. and Pittet, C., Isoperimetry and Volume Growth for Cayley Graphs. Unfinished manuscript, (2014). Available at: https://pi.math.cornell.edu/~lsc/papers/surv.pdf}.  
	
	\item {\small Tedeev, A.F., The interface blow-up phenomenon and local estimates for doubly degenerate parabolic equations. Appl. Anal. \textbf{86}(6) (2007), 755–782. DOI: 10.1080/00036810701435711.}
	
	\item {\small Ting, G.Y. and Feng, W.L., p-Laplace elliptic inequalities on the graph. Commun. Pure Appl. Anal. \textbf{8}(3) (2025), 389–411. DOI: 10.3934/cpaa.2024094.}
	
	\item {\small Xin, Q., Mu, C., and Liu, D., Extinction and Positivity of the Solutions for a p-Laplacian Equation with Absorption on Graphs. J. Appl. Math. (2011), Art. ID 937079. DOI: 10.1155/2011/937079.}
\end{enumerate}

{\small \bigskip}

{\small \bigskip}

{\small \bigskip}

{\small \bigskip}

\end{document}